\documentclass[11pt,reqno]{amsart}
\usepackage[T1]{fontenc}
\usepackage[utf8]{inputenc}
\usepackage[english]{babel}
\usepackage{lmodern,microtype}
\usepackage[a4paper,margin=29mm]{geometry}
\usepackage{amsmath,amssymb,amsthm,mathtools,booktabs,enumitem}
\usepackage{xcolor,hyperref,needspace}
\definecolor{linkblue}{RGB}{30,65,100}
\hypersetup{colorlinks=true,linkcolor=linkblue,citecolor=linkblue,urlcolor=linkblue,
 pdftitle={Cubic fourfolds with a symplectic involution from finite Prym data},pdfauthor={Lucas Li Bassi}}
\newtheorem{theorem}{Theorem}[section]
\newtheorem{proposition}[theorem]{Proposition}
\newtheorem{lemma}[theorem]{Lemma}

\newtheorem{maintheorem}{Theorem}
\renewcommand{\themaintheorem}{\Alph{maintheorem}}
\newtheorem{maincorollary}{Corollary}[maintheorem]
\newtheorem*{corollaryA}{Corollary A.1}
\newtheorem*{theoremA}{Theorem A}
\newtheorem*{theoremB}{Theorem B}
\theoremstyle{definition}

\newtheorem{example}[theorem]{Example}
\theoremstyle{remark}
\newtheorem{remark}[theorem]{Remark}
\newcommand{\CC}{\mathbb C}
\newcommand{\ZZ}{\mathbb Z}
\newcommand{\QQ}{\mathbb Q}
\newcommand{\PP}{\mathbb P}

\DeclareMathOperator{\Aut}{Aut}

\DeclareMathOperator{\NS}{NS}
\DeclareMathOperator{\im}{im}
\DeclareMathOperator{\Sym}{Sym}
\DeclareMathOperator{\rk}{rk}
\DeclareMathOperator{\Fix}{Fix}
\numberwithin{equation}{section}
\newcommand{\sectionintrospace}{\par\addvspace{2\baselineskip}}

\title[Cubic fourfolds and Prym data]{Cubic fourfolds with a symplectic involution from finite Prym data}
\author{Lucas Li Bassi}
\date{\today}
\address{ICERM at Brown University, 121 Main Street, Providence, RI 02903, USA}
\email{lucas.libassi@gmail.com}

\begin{document}
\begin{abstract}
We reconstruct smooth cubic fourfolds with a symplectic involution through classical interpolation of smooth plane conics from five distinct points, or four with their cross-ratio. If the quadratic coefficient map has rank three and the fixed cubic surface has trivial automorphism group, five polarized anti-invariant Jacobian factors of smooth Eckardt sections suffice, or four with their parameters. Both bounds are optimal. By contrast, in rank two we construct two-dimensional families of equivariantly nonisomorphic fourfolds with the same entire unparametrized Prym image. The factor family also has an integral projective Lagrangian closure in an associated OG10 compactification.
\end{abstract}
\maketitle
\tableofcontents

\section{Introduction}

The intermediate Jacobians of the smooth hyperplane sections of a cubic fourfold form a family over an open subset of the dual projective space. After passing to an associated OG10 compactification, this family determines the cubic. Indeed, the resulting Lagrangian fibration recovers its primitive Hodge lattice \cite[Theorem~4.6]{DMS}, and the Torelli theorem for cubic fourfolds then gives its projective isomorphism class. This raises the question of how much of the Jacobian family is needed for reconstruction.

In this paper, we study this question for cubic fourfolds with a symplectic involution. Indeed, the involution supplies a distinguished pencil of hyperplane sections, and we can consider the polarized anti-invariant factors of their intermediate Jacobians. It is then natural to ask whether finitely many such factors determine the fourfold together with its involution, and whether knowing their positions in the pencil reduces the number required.

To describe these data, let $V$ and $W$ be complex vector spaces of dimensions four and two. The equation of a smooth cubic fourfold with a symplectic involution can be written
\begin{equation}\label{eq:normalform}
Y=\{F(x)+q(x,z)=0\}\subset\PP(V\oplus W),
\qquad \iota[x:z]=[x:-z],
\end{equation}
where $F\in\Sym^3V^*$ and $q\in V^*\otimes\Sym^2W^*$. Its fixed locus consists of the smooth cubic surface $S=\{F=0\}\subset\PP(V)$ and the line $\PP(W)$ \cite[Section~4.2]{Mar}. The hyperplanes containing $\PP(V)$ therefore give a pencil of invariant sections,
\[
T_t=Y\cap\PP(V\oplus\CC t),\qquad [t]\in\PP(W),
\]
called the \emph{Eckardt pencil}. For each smooth member, the induced involution $\iota_t=\iota|_{T_t}$ defines an anti-invariant factor
\[
A_t=\im(1-\iota_t:J(T_t)\longrightarrow J(T_t)),
\qquad \lambda_t=\theta_t|_{A_t},
\]
where $\theta_t$ is the principal polarization of $J(T_t)$. Thus our data are polarized abelian fourfolds of type $(1,1,1,2)$ \cite[Lemma~1.12]{CMZ}, considered up to polarized isomorphism.

To relate the factors to the equation of $Y$, we use the Torelli theorem of Casalaina-Martin--Zhang, recalled in Theorem~\ref{thm:cmz}: each $(A_t,\lambda_t)$ determines the surface--plane pair $(S,\{q(-,t)=0\})$ \cite[Theorem~6.1]{CMZ}. If $\Aut(S)=1$, the recovered copies of $S$ can be identified uniquely, and their anticanonical embeddings place all the planes in the same dual projective space. We can therefore place the recovered planes on the image of
\[
c_q:\PP(W)\longrightarrow\PP(V^*),\qquad [t]\longmapsto[q(-,t)].
\]
This morphism is defined everywhere, since smoothness of $Y$ forces $q(-,t)\ne0$ for $t\ne0$.

Let $\widehat q:\Sym^2W\to V^*$ be the quadratic coefficient map. When $\rk\widehat q=3$, the morphism $c_q$ embeds the parameter line as a smooth conic $C_q$. The surface determines $F$ up to scalar, while a parametrization of this conic determines $q$ up to a change of the $W$ coordinates. Thus the pair $(S,C_q)$ determines the equivariant fourfold (Proposition~\ref{prop:surfaceconic}), and finite reconstruction reduces to classical interpolation for plane conics. Five factors suffice because five distinct points determine the conic. Four points instead leave a pencil of conics, but knowing the corresponding parameters specifies their ordered cross-ratio and singles out one smooth member. This gives the following two forms of finite reconstruction, where we use primed notation for a second fourfold.

\begin{maintheorem}\label{thm:reconstruction}
Let $(Y,\iota)$ and $(Y',\iota')$ be smooth cubic fourfolds of the form~\eqref{eq:normalform}, with $\rk\widehat q=\rk\widehat q'=3$ and $\Aut(S)=\Aut(S')=1$. Then there is a projective isomorphism $f:Y\xrightarrow{\sim}Y'$ satisfying $f\circ\iota=\iota'\circ f$ in either of the following cases.
\begin{enumerate}[label=\textup{(\roman*)},ref=\themaintheorem(\roman*)]
\item\label{thm:five} Five distinct smooth sections in each Eckardt pencil give the same unordered set of polarized isomorphism classes of anti-invariant factors.
\item\label{thm:four} There are four distinct smooth Eckardt parameters $t_1,\ldots,t_4$ for $Y$ and $t'_1,\ldots,t'_4$ for $Y'$, and a projectivity of the parameter lines sending $t_i$ to $t'_i$, such that
\[
(A_{t_i},\lambda_{t_i})\simeq(A'_{t'_i},\lambda'_{t'_i})
\qquad (i=1,\ldots,4)
\]
as polarized abelian varieties.
\end{enumerate}
\end{maintheorem}

In part~\textup{(ii)}, each factor is attached to its parameter, and the parameter configuration is considered up to a common projectivity. Our ordered cross-ratio is the image of the fourth point after sending the first three to $0,1,\infty$.

\begin{maincorollary}\label{cor:optimal}\label{cor:five-optimal}\label{cor:four-optimal}
Both bounds in Theorem~\ref{thm:reconstruction} are optimal. There are one-dimensional families of pairwise equivariantly nonisomorphic smooth rank-three fourfolds with a fixed surface $S$ satisfying $\Aut(S)=1$ and four polarized factors in common. There are also such families with the same three polarized factors at three fixed parameters.
\end{maincorollary}

The rank-three assumption is essential to this argument. In rank two, $c_q$ is a double cover of a line, so its image no longer records its branch pair. We can therefore keep the surface and image line fixed while varying the covering. This gives families for which even the whole unparametrized set of factors fails to determine the fourfold.

\begin{maintheorem}\label{thm:ranktwo}
There is a two-dimensional algebraic family of pairwise equivariantly nonisomorphic smooth cubic fourfolds with symplectic involutions. Its members have $\rk\widehat q=2$, the same fixed cubic surface $S$ with $\Aut(S)=1$, and the same entire set of polarized isomorphism classes of anti-invariant factors from smooth Eckardt sections.
\end{maintheorem}

After proving the reconstruction results in Sections~\ref{sec:geometry}--\ref{sec:ranktwo}, we treat torsors and surfaces with automorphisms in Section~\ref{sec:extensions}. Section~\ref{sec:lsv} studies the full Jacobian fibration and the integral projective Lagrangian closure of the factor family in an associated OG10 compactification (Proposition~\ref{prop:prymclosure}).

\section*{Acknowledgements}
I thank Chiara Camere for useful conversations on this topic. Moreover I would like to thank Hannah Dell for having spotted some inaccuracies on a previous draft and Alessandra Sarti for her suggestions on how to improve the exposition. This material is based upon work supported by the National Science Foundation under Grant No. DMS-2424556 while the author was in residence at the Institute for Computational and Experimental Research in Mathematics in Providence, RI, during the Computations on K3 Surfaces and Related Varieties semester program.
\subsection*{AI disclosure}
I used an AI assistant based on a large language model (LLM) to identify relevant references, improve the exposition and assist with computations. In particular, it found the explicit expression in Equation~\eqref{eq: nonzero pol}. All other mathematical ideas were developed independently of the use of AI. I have independently verified the mathematical content and take full responsibility for any remaining errors.

\Needspace{12\baselineskip}
\section{The cubic surface and its conic}
\label{sec:geometry}

The fixed surface determines the cubic term in~\eqref{eq:normalform}, so reconstruction amounts to recovering the quadratic term. We describe it through the planes associated with the Eckardt sections and show that, in rank three, their image conic determines it up to a change of coordinates. We also establish the smoothness conditions needed to vary these data in the optimality argument.
\sectionintrospace

We work over $\CC$. Every isomorphism of smooth cubic fourfolds extends to a projectivity: indeed, $\operatorname{Pic}(Y)=\ZZ[\mathcal O_Y(1)]$ and $-K_Y=3\mathcal O_Y(1)$, so the hyperplane class is preserved. Likewise, an isomorphism of smooth cubic surfaces preserves their anticanonical embeddings. We may therefore compare both varieties through their equations.

\subsection{The invariant equation}
Consider the involution $\iota[x:z]=[x:-z]$ on $\PP(V\oplus W)$, where $\dim V=4$ and $\dim W=2$. If it preserves a smooth cubic fourfold, the defining polynomial $P$ satisfies $P(x,-z)=\pm P(x,z)$. Smoothness excludes the minus sign: in that case the equation would have the form
\[
z_0Q_0(x)+z_1Q_1(x)+G(z_0,z_1).
\]
Such a fourfold would contain $\PP(V)$ and be singular at the nonempty common zero locus of $Q_0,Q_1$ there. We therefore obtain the invariant equation~\eqref{eq:normalform}; see also \cite[Section~7]{Camere}, \cite[Section~4.2]{Mar} and \cite[Section~3.1]{FFM}. Moreover, the lift $\mathrm{id}_V\oplus(-\mathrm{id}_W)$ has determinant $1$ and preserves this polynomial. It follows that it acts trivially on $H^{3,1}(Y)$, so the induced involution on $Y$ is symplectic.

Smoothness also describes the fixed locus. At $[x:0]\in Y$, the gradient is $(dF_x,0)$, so the surface $S=\{F=0\}\subset\PP(V)$ is smooth. At $[0:z]\in\PP(W)$, it is $(q(-,z),0)$, and hence $q(-,z)\ne0$ for $z\ne0$. We therefore obtain
\begin{equation}\label{eq:fixedsystems}
\Fix(\iota|_Y)=S\sqcup\PP(W),\qquad
(\PP(V\oplus W)^*)^\iota=\PP(V^*)\sqcup\PP(W^*).
\end{equation}
The two components in the dual space are the invariant hyperplane systems of \cite[Section~4.5]{MarTh}. We use the component $\PP(W^*)$, which parametrizes the Eckardt pencil and whose sections all contain $S$. By contrast, a smooth section from the other component carries a non-Eckardt involution. Its invariant Jacobian factor has dimension three and polarization type $(1,2,2)$, and the Torelli theorem of Casalaina-Martin--Marquand--Zhang reconstructs the equivariant threefold from this factor \cite[Lemma~1.7 and Theorem~3.1]{CMMZ}.

To describe the Eckardt sections, choose $t\in W\setminus\{0\}$ and set $z=ut$. Then
\begin{equation}\label{eq:eckardtsection}
T_t=\{F(x)+u^2\ell_t(x)=0\},\qquad \ell_t=q(-,t).
\end{equation}
Rescaling $t$ only rescales $\ell_t$, so the plane $\{\ell_t=0\}$ depends only on $[t]$. When $T_t$ is smooth, $p_t=[0:1]$ is an Eckardt point: its tangent hyperplane section is a cone over the smooth plane cubic $S\cap\{\ell_t=0\}$, with vertex $p_t$ \cite[Definition~1.1 and Proposition~1.2]{CMZ}. The smoothness of this plane cubic also characterizes smoothness of $T_t$.

\begin{lemma}\label{lem:smoothsection}
Let $S=\{F=0\}$ be a smooth cubic surface and let $\ell$ be a nonzero linear form. Then the cubic threefold $T=\{F+u^2\ell=0\}$ is smooth if and only if the plane $\{\ell=0\}$ meets $S$ transversely.
\end{lemma}
\begin{proof}
Since $\ell\ne0$, the gradient in the $x$ directions does not vanish at $[0:1]$. The threefold is also smooth along $u=0$, by smoothness of $S$. We are therefore left to consider points with $x\ne0$ and $u\ne0$. At such a point, the equation and its derivatives give the singularity conditions
\[
F(x)=\ell(x)=0,\qquad dF_x+u^2\ell=0.
\]
The last equality says that $\{\ell=0\}$ is tangent to $S$ at $[x]$. Conversely, if this plane is tangent to $S$ at $[x]$, then $dF_x$ is a nonzero scalar multiple of $\ell$. Over $\CC$ we may choose $u$ with the required square, so all the displayed conditions hold and $T$ has a singular point.
\end{proof}

We denote the open set of transverse planes by
\[
\Omega_S=\{[\ell]\in\PP(V^*):\{F+u^2\ell=0\}\text{ is smooth}\}
         =\PP(V^*)\setminus S^\vee.
\]
Lemma~\ref{lem:smoothsection} shows that this condition depends only on $[\ell]$. Thus, whenever $q$ has no base points, the map $c_q([t])=[q(-,t)]$ is defined on $\PP(W)$, and the set $U_q$ of smooth Eckardt parameters satisfies
\begin{equation}\label{eq:smoothimage}
U_q=c_q^{-1}(\Omega_S),\qquad
c_q(U_q)=c_q(\PP(W))\cap\Omega_S.
\end{equation}
The second equality records the planes which occur, but forgets how many parameters map to each plane and whether it is a branch value of $c_q$.

\subsection{Recovering the equation from a conic}
Suppose now that $\rk\widehat q=3$. Since $\dim\Sym^2W=3$, the coefficient map is injective. Consequently,
\begin{equation}\label{eq:conicmap}
c_q:\PP(W)\longrightarrow\PP(V^*),\qquad [t]\longmapsto[\ell_t]
\end{equation}
is the Veronese embedding followed by a linear embedding, and its image $C_q$ is a smooth plane conic.

A parametrization of $C_q$ gives the quadratic coordinates of $q$. We must therefore check that changing the parametrization or its common scalar does not change the equivariant isomorphism class. In the pair $(S,C_q)$, we retain the position of the conic in the dual of the ambient space of $S$.

\begin{proposition}\label{prop:surfaceconic}
Let $(Y,\iota)$ and $(Y',\iota')$ be smooth cubic fourfolds of the form~\eqref{eq:normalform}, with $\rk\widehat q=\rk\widehat q'=3$. There is a projective equivariant isomorphism $Y\simeq Y'$ if and only if there is a projectivity
\[
a:\PP(V)\xrightarrow{\sim}\PP(V')
\]
which sends $S$ to $S'$ and whose induced map on hyperplanes, $H\mapsto a(H)$, sends $C_q\subset\PP(V^*)$ to $C_{q'}\subset\PP(V'^*)$.
\end{proposition}
\begin{proof}
Let $\varphi:Y\to Y'$ be an equivariant isomorphism. Since it preserves the hyperplane class, it extends to a projectivity
\[
\Phi:\PP(V\oplus W)\xrightarrow{\sim}\PP(V'\oplus W').
\]
To see that $\Phi$ is equivariant on the ambient spaces, choose matrices $A,B$ representing $\Phi\circ\iota$ and $\iota'\circ\Phi$. These projectivities agree on $Y$, and therefore the quadrics
\[
(Ay)_i(By)_j-(Ay)_j(By)_i,
\qquad 0\leq i<j\leq5,
\]
vanish on $Y$. However, the ideal of this reduced cubic hypersurface is generated by its cubic equation and therefore contains no nonzero quadric. The displayed polynomials must vanish identically, so $Ay$ and $By$ are proportional for every nonzero $y$. This proves that
\[
\Phi\circ\iota=\iota'\circ\Phi
\]
on the ambient projective spaces.

It follows that a linear lift of $\Phi$ conjugates the involutions up to scalar. Since their eigenspaces have different dimensions, four and two, it cannot exchange them. The lift thus gives a projectivity $a$ sending $S$ to $S'$, together with a change of the $W$ coordinates. The latter reparametrizes the Eckardt pencil, and the induced action of $a$ on its associated planes sends $C_q$ to $C_{q'}$.

Conversely, suppose that the projectivity $a$ is given. After lifting it to a linear map and identifying $V$ with $V'$, we have a single surface--conic pair $(S,C)$. The two equations defining $S$ differ by a nonzero scalar, so we may also normalize the fourfold equations to have the same cubic term $F$. It remains to compare their quadratic terms.

To do so, choose an isomorphism $\PP^1\to C$. Since the pullback of $\mathcal O_C(1)$ is $\mathcal O_{\PP^1}(2)$, this parametrization can be written as
\[
[u:v]\longmapsto[u^2\ell_0+uv\ell_1+v^2\ell_2],
\]
with independent $\ell_0,\ell_1,\ell_2\in V^*$. In particular, it gives a quadratic term
\[
q(x,u,v)=u^2\ell_0(x)+uv\ell_1(x)+v^2\ell_2(x).
\]
Now $c_q$ and $c_{q'}$ are two parametrizations of the same conic, so they differ by an element of $\mathrm{PGL}_2$, which lifts to a change of the $W$ coordinates. After making this change, their base-point-free tuples of quadratic forms define the same morphism and hence the same pullback of the hyperplane bundle. The two identifications of this pullback with $\mathcal O_{\PP^1}(2)$ differ by a global unit on $\PP^1$. Therefore $q'=\mu q$ for some $\mu\in\CC^*$, and rescaling the $W$ coordinates by $\mu^{-1/2}$ identifies $\{F+q=0\}$ with $\{F+q'=0\}$. Since each coordinate change preserves the eigenspaces, the resulting isomorphism is equivariant.
\end{proof}

The isomorphism in the proposition need not be unique: both the identity of $Y$ and $\iota$ act trivially on $S$ and $C_q$.

\begin{remark}
This description is dual to Marquand's description by a cubic surface and a quadric cone \cite[Proposition~3.5.1]{MarTh}. Indeed, if $q=\ell_0z_0^2+\ell_1z_0z_1+\ell_2z_1^2$, the cone has equation $\ell_1^2-4\ell_0\ell_2=0$ in $\PP(V)$, and its dual conic is precisely $C_q$. The argument above uses the equation directly, so it also applies when the branch curve on the cone is not smooth.
\end{remark}

To construct the optimality families, we will vary the conic while keeping the surface fixed. The following lemma ensures that we can start from a smooth fourfold for every smooth cubic surface.

\begin{lemma}\label{lem:nonempty}
Let $\dim V=4$ and $\dim W=2$.
\begin{enumerate}[label=\textup{(\roman*)}]
\item For each smooth cubic surface $S=\{F=0\}\subset\PP(V)$, the terms $q\in V^*\otimes\Sym^2W^*$ for which
\[
Y_q=\{F(x)+q(x,z)=0\}\subset\PP(V\oplus W)
\]
is smooth and $\rk\widehat q=3$ form a nonempty Zariski-open subset.
\item There is a nonempty Zariski-open subset of the moduli space of smooth cubic surfaces on which $\Aut(S)=1$.
\end{enumerate}
In particular, there are smooth fourfolds of the form~\eqref{eq:normalform} with $\rk\widehat q=3$ and $\Aut(S)=1$.
\end{lemma}
\begin{proof}
For (i), fix $F$ and set $Q=V^*\otimes\Sym^2W^*$. The linear system
\[
\PP(\CC F\oplus Q)
\]
has equations $aF(x)+q(x,z)=0$ and set-theoretic base locus
\[
S\sqcup\PP(W).
\]
Indeed, the monomials $x_i z_jz_k$ vanish simultaneously precisely when $x=0$ or $z=0$; on the latter locus, the equation $F=0$ cuts out $S$.

Bertini's theorem gives smoothness of a general member away from this base locus. Since $a\ne0$ and $\rk\widehat q=3$ are nonempty open conditions on the irreducible parameter space $\PP(\CC F\oplus Q)$, we may impose both conditions as well. They allow us to check smoothness along the base locus: at $[x:0]\in S$, the gradient is
\[
(a\,dF_x,0),
\]
which is nonzero because $a\ne0$ and $S$ is smooth. At $[0:z]\in\PP(W)$, the gradient is instead
\[
(q(-,z),0).
\]
Since $\widehat q$ is injective, we have $q(-,z)=\widehat q(z^2)\ne0$ for $z\ne0$. The member is therefore smooth everywhere. Dividing its equation by $a$, we obtain the required form $F+q/a=0$ and hence a point of the set in (i). This set is open because smoothness and maximal rank are open conditions on $Q$.

For (ii), we use the classification of Dolgachev--Duncan \cite[Theorem~1.1 and Table~1]{DD}, which holds over algebraically closed fields of arbitrary characteristic. Over $\CC$, the surfaces with nontrivial automorphisms lie in finitely many strata of dimension at most three. Since the moduli space is irreducible of dimension four, the closures of these strata are proper closed subsets. Their complement is consequently a nonempty open set on which $\Aut(S)=1$. Choosing $S$ in this complement and applying (i) proves the last assertion.
\end{proof}

We also need sufficiently many smooth members of the Eckardt pencil. Bertini's theorem gives smoothness of a general member away from its base surface $S$. Along $S$, the differential of the hyperplane equation involves only the $z$ variables, whereas $dF$ involves only the $x$ variables. These differentials are therefore independent, so the section is smooth along $S$ as well. Thus the smooth parameters form a nonempty open set; equivalently, $C_q\not\subset S^\vee$. In particular, we can choose as many distinct smooth parameters as the reconstruction statements require.

\begin{example}
One can also see the nonemptiness of the smooth rank-three locus from the equation
\begin{equation}\label{eq:example}
F=x_0^3+x_1^3+x_2^3+x_3^3,
\qquad q=x_0z_0^2+x_1z_0z_1+x_2z_1^2.
\end{equation}
To check smoothness, suppose that all derivatives vanish at a point. The derivative in $x_3$ first gives $x_3=0$. If either $z$ coordinate were zero, the other derivative equations would force all coordinates to vanish, which is impossible. We may therefore set $z_0=1$ and $z_1=t\ne0$. The derivatives in the $z$ variables now give $x_1=-2x_0/t$ and $x_2=x_0/t^2$. Substituting into the remaining equations yields $3x_0^2=-1$, $t^3=4$ and $t^6=1$. The last two equalities are incompatible, proving that the fourfold is smooth.

The fixed surface in this example is the Fermat surface, so its automorphism group is nontrivial. We can nevertheless obtain the remaining hypothesis by keeping $q$ fixed and deforming $F$: smoothness is open, while a general deformation lies in the open locus $\Aut(S)=1$ of Lemma~\ref{lem:nonempty}.

\end{example}

\Needspace{12\baselineskip}
\section{Finite reconstruction and parameter data}
\label{sec:five}

By Proposition~\ref{prop:surfaceconic}, we need to recover the fixed surface and its embedded conic. The Torelli theorem for cubic surface pairs recovers the surface together with one plane from each factor. We will use $\Aut(S)=1$ to identify these copies of the surface, and then apply conic interpolation to the recovered planes. This proves both parts of Theorem~\ref{thm:reconstruction}; the remaining pencil of conics supplies the optimality family.
\sectionintrospace

For a smooth Eckardt section $T_t$, let $\lambda_t$ be the restriction of the principal polarization to $A_t=\im(1-\iota_t)$. We use the following form of the Torelli theorem of Casalaina-Martin--Zhang.

\Needspace{8\baselineskip}
\begin{theorem}[Casalaina-Martin--Zhang]\label{thm:cmz}
For a smooth section $T_t=\{F+u^2\ell_t=0\}$, the polarized abelian variety $(A_t,\lambda_t)$ has dimension four and polarization type $(1,1,1,2)$. Its polarized isomorphism class determines the pair $(S,\{\ell_t=0\})$ up to projective isomorphism.
\end{theorem}
\begin{proof}
The dimension and polarization type are established in \cite[Lemma~1.12]{CMZ}. For the reconstruction assertion, apply the global injectivity of \cite[Theorem~6.1]{CMZ}, which holds for every smooth cubic surface with a transverse plane section. By Lemma~\ref{lem:smoothsection}, these hypotheses follow from the smoothness of $T_t$.
\end{proof}

Since the recovered pairs are defined up to projective isomorphism, we use $\Aut(S)=1$ to identify their surfaces.

\begin{lemma}\label{lem:alignment}
Assume that $\Aut(S)=1$. The polarized isomorphism classes of a finite collection of factors determine their associated planes as points of
\[
\PP\bigl(H^0(S,-K_S)\bigr),
\]
up to a common projective isomorphism.
\end{lemma}
\begin{proof}
Write $(S_i,\Pi_i)$ for the pairs given by Theorem~\ref{thm:cmz}, and choose a representative $S$ of their common surface isomorphism class. Since any two isomorphisms $S_i\to S$ differ by an automorphism of $S$, the assumption $\Aut(S)=1$ makes each such identification unique.

These identifications also transport the planes, because the embeddings of the surfaces are anticanonical. Indeed, adjunction gives $-K_S\simeq\mathcal O_S(1)$, and restriction identifies $V^*$ with $H^0(S,\mathcal O_S(1))$. Thus $\PP(H^0(S,-K_S))\simeq\PP(V^*)$, where $[\ell|_S]$ represents the plane $\{\ell=0\}$. It follows that all the recovered planes are placed in this single dual space. Replacing the representative $S$ transports them by one common projectivity, which is exactly the ambiguity in the statement.
\end{proof}

\subsection{Conics through prescribed points}

The recovered planes are now points in a common dual space, so reconstruction reduces to the following classical facts about plane conics. We include a proof for the reader's convenience; the explicit parametrization will also be used in the optimality family.

\Needspace{5\baselineskip}
\begin{lemma}\label{lem:conicinterpolation}
Let $C$ be a smooth plane conic in a projective space.
\begin{enumerate}[label=\textup{(\roman*)}]
\item Five distinct points on $C$ determine it uniquely.
\item Four distinct points on $C$ determine a pencil of conics in their plane. Among its smooth members, the ordered cross-ratio of the four points determines the conic uniquely.
\end{enumerate}
\end{lemma}
\begin{proof}
Since three distinct points of $C$ are not collinear, the prescribed points already determine its plane. Any other conic through them therefore lies in that plane, and its equation restricts to a section of $\mathcal O_C(2)\simeq\mathcal O_{\PP^1}(4)$. If there are five distinct common points, the restriction has five zeros and must vanish identically. The two conics are consequently equal, proving (i).

For (ii), choose coordinates in which the four ordered points are $e_0,e_1,e_2,(1:1:1)$. Passing through the first three removes the square terms from the quadratic equation, while passing through the fourth imposes one relation on the remaining coefficients. Thus the conics through the four points form the pencil
\begin{equation}\label{eq:conicpencil}
aXY+bXZ+cYZ=0,\qquad a+b+c=0.
\end{equation}
A member is smooth precisely when $abc\ne0$, so we can divide the equation by the coefficient of $XY$. The smooth members then take the form
\[
C_s=\{XY+sXZ-(1+s)YZ=0\},\qquad s\in\CC\setminus\{0,-1\}.
\]
To compute the cross-ratio on $C_s$, use the parametrization
\begin{equation}\label{eq:conicparam}
c_s([u:v])=[(1+s)uv:u^2+suv:uv+sv^2].
\end{equation}
The coefficient matrix has determinant $-s(1+s)\ne0$, so $c_s$ is an isomorphism onto $C_s$. The four marked points have parameters $-s,\infty,0,1$ in the coordinate $t=u/v$. The projectivity
\[
t\longmapsto\frac{t+s}{t}
\]
sends the first three to $0,1,\infty$ and the fourth to $1+s$. Therefore their ordered cross-ratio is $r=1+s$. It determines $s$, and hence the conic, as required.
\end{proof}

\subsection{Reconstruction and optimality}
\label{sec:sharpness}

We now apply Lemma~\ref{lem:conicinterpolation} to the planes recovered by the Torelli theorem.

\begin{theoremA}
Let $(Y,\iota)$ and $(Y',\iota')$ be smooth cubic fourfolds of the form~\eqref{eq:normalform}, with $\rk\widehat q=\rk\widehat q'=3$ and $\Aut(S)=\Aut(S')=1$. Then there is a projective isomorphism $f:Y\xrightarrow{\sim}Y'$ satisfying $f\circ\iota=\iota'\circ f$ in either of the following cases.
\begin{enumerate}[label=\textup{(\roman*)}]
\item Five distinct smooth sections in each Eckardt pencil give the same unordered set of polarized isomorphism classes of anti-invariant factors.
\item There are four distinct smooth Eckardt parameters $t_1,\ldots,t_4$ for $Y$ and $t'_1,\ldots,t'_4$ for $Y'$, and a projectivity of the parameter lines sending $t_i$ to $t'_i$, such that
\[
(A_{t_i},\lambda_{t_i})\simeq(A'_{t'_i},\lambda'_{t'_i})
\qquad (i=1,\ldots,4)
\]
as polarized abelian varieties.
\end{enumerate}
\end{theoremA}
\begin{proof}
By Theorem~\ref{thm:cmz}, corresponding factors determine the same surface--plane pairs. Lemma~\ref{lem:alignment} identifies their surfaces and places the planes in a common dual space. These planes are distinct: indeed, $c_q$ is an embedding, and an isomorphism between two pairs on $S$ would be an automorphism carrying one plane to the other. Since $\Aut(S)=1$, distinct parameters therefore give distinct pair classes.

In case~\textup{(i)}, the two image conics contain the same five points, so they coincide by Lemma~\ref{lem:conicinterpolation}(i).

In case~\textup{(ii)}, use the given projectivity to identify the parameter lines. The four corresponding factors then give the same ordered points $p_1,\ldots,p_4$. Since $c_q$ and $c_{q'}$ are isomorphisms onto their conics, each transports the cross-ratio of the specified parameters to the cross-ratio of these four points. The conics consequently contain the same points with the same ordered cross-ratio, and they coincide by Lemma~\ref{lem:conicinterpolation}(ii).

In either case, the two fourfolds have the same surface--conic pair. Proposition~\ref{prop:surfaceconic} therefore gives the required equivariant projective isomorphism.
\end{proof}

To prove optimality, we keep four surface--plane pairs fixed and vary the conic through their four points. The parametrization in~\eqref{eq:conicparam} will also preserve three of their parameters, so the same construction treats both bounds.

\begin{corollaryA}
Both bounds in Theorem~\ref{thm:reconstruction} are optimal. There are one-dimensional families of pairwise equivariantly nonisomorphic smooth rank-three fourfolds with a fixed surface $S$ satisfying $\Aut(S)=1$ and four polarized factors in common. There are also such families with the same three polarized factors at three fixed parameters.
\end{corollaryA}
\begin{proof}
Choose a smooth rank-three fourfold with $\Aut(S)=1$ as in Lemma~\ref{lem:nonempty}, and write $(S,C)$ for its surface--conic pair. Fix four points of $C$ corresponding to smooth Eckardt sections. In the coordinates of~\eqref{eq:conicparam}, let $\ell_X,\ell_Y,\ell_Z$ be the corresponding basis of the three-dimensional subspace of $V^*$. Set
\begin{equation}\label{eq:sharpfamily}
q_s=(1+s)uv\,\ell_X+(u^2+suv)\ell_Y+(uv+sv^2)\ell_Z.
\end{equation}
For $s\ne0,-1$, these terms give an algebraic family $Y_s=\{F+q_s=0\}$ of rank-three cubics. If $C=C_{s_0}$, the original parametrization of $C$ differs from $c_{s_0}$ by a projectivity of $\PP(W)$. Lifting it and rescaling the $W$ coordinates as in Proposition~\ref{prop:surfaceconic} identifies $Y_{s_0}$ with the chosen smooth fourfold. Smoothness is open, so the family is smooth over a nonempty Zariski-open set containing $s_0$.

The four fixed points give the same surface--plane pairs throughout this open set, and hence the same four polarized factors. Nevertheless, the fourfolds are pairwise equivariantly nonisomorphic. Indeed, an equivariant isomorphism would restrict to the identity on $S$, because $\Aut(S)=1$, and would therefore act trivially on $\PP(V)$ and its dual. Proposition~\ref{prop:surfaceconic} would then identify the conics, whereas distinct values of $s$ give distinct members of the pencil. This proves optimality in part~\textup{(i)}.

For part~\textup{(ii)}, observe that the parameters $[1:0]$, $[0:1]$ and $[1:1]$ have the fixed images $e_1$, $e_2$ and $(1:1:1)$ throughout the family. These images correspond to smooth plane sections of $S$. On the same open set, we therefore have pairwise equivariantly nonisomorphic smooth fourfolds with three fixed parameters and the same three polarized factors.
\end{proof}

\begin{remark}
The theorem uses only polarized isomorphism classes, without embeddings in intermediate Jacobians or a common marking.
\end{remark}

\begin{remark}
One can also check the dimensions of the two Jacobian factors directly from the topological Lefschetz formula. Indeed, the fixed locus of $\iota_t$ on $T_t$ is $S\sqcup\{p_t\}$ and has Euler characteristic $10$, while the even cohomology contributes $4$. The formula therefore gives
\[
\operatorname{tr}(\iota_t|H^3(T_t,\QQ))=-6.
\]
Combining this trace with $b_3(T_t)=10$, we obtain invariant and anti-invariant ranks $2$ and $8$, and hence abelian dimensions $1$ and $4$. Thus the fourfold $A_t$ is the anti-invariant factor, in agreement with Theorem~\ref{thm:cmz}, whereas the identity component of the fixed locus in the Jacobian has dimension one.

The same calculation applies to $Y$. Its fixed locus has Euler characteristic $9+2=11$, and the four cohomology groups outside the middle degree contribute $4$. After removing the hyperplane-square class, the trace on primitive $H^4$ is therefore $6$, giving invariant and anti-invariant ranks $14$ and $8$.
\end{remark}

\Needspace{12\baselineskip}
\section{Rank two and the branch values}
\label{sec:ranktwo}

In rank two, the parameter map is a double cover of a line, so its image no longer determines its branch pair. We can therefore vary the covering while keeping the surface--plane pairs unchanged. We now construct this variation algebraically and show that its members are the nonisomorphic fourfolds of Theorem~\ref{thm:ranktwo}.
\sectionintrospace

Fix a smooth cubic surface $S=\{F=0\}\subset\PP^3$ with $\Aut(S)=1$. By Bertini's theorem, we can choose a line $\Lambda\subset(\PP^3)^\vee$ whose pencil of planes has base line
\[
B=\bigcap_{H\in\Lambda}H
\]
transverse to $S$. The same theorem ensures that a general plane of this pencil meets $S$ transversely away from $B$. Moreover, if a plane of the pencil were tangent to $S$ at a point of $B$, the tangent plane to $S$ there would contain $B$, contrary to transversality. Thus a general plane of the pencil is transverse to $S$ everywhere, and the set $\Delta=\Lambda\cap S^\vee$ of nontransverse planes is finite.

Given distinct planes $H_0,H_1\in\Lambda\setminus\Delta$, with equations $\ell_0,\ell_1$, set
\begin{equation}\label{eq:ranktwo}
Y_{H_0,H_1}=\{F(x)+\ell_0(x)u^2+\ell_1(x)v^2=0\},
\qquad \iota[x:u:v]=[x:-u:-v].
\end{equation}
This construction depends, up to equivariant isomorphism, only on the unordered pair $\{H_0,H_1\}$. Indeed, rescaling either plane equation can be absorbed by rescaling the corresponding coordinate $u$ or $v$, while interchanging the planes amounts to exchanging these coordinates. To use these pairs as parameters for smooth fourfolds, we must also check that the transversality conditions exclude singularities.

\begin{lemma}\label{lem:ranktwosmooth}
Let $S$ be a smooth cubic surface and let $H_0,H_1$ be distinct planes, with independent equations $\ell_0,\ell_1$. If $H_0$, $H_1$ and their common line $B$ all meet $S$ transversely, then the fourfold in~\eqref{eq:ranktwo} is smooth.
\end{lemma}
\begin{proof}
The fourfold is smooth along $u=v=0$ because $S$ is smooth. It is also smooth at $x=0$: the gradient in the $x$ directions is $u^2\ell_0+v^2\ell_1$, which is nonzero by the independence of $\ell_0$ and $\ell_1$.

It remains to consider $x\ne0$ with at least one of $u,v$ nonzero. The derivatives in these variables are $2u\ell_0(x)$ and $2v\ell_1(x)$. If exactly one variable is nonzero, their vanishing at a singular point would make the corresponding plane tangent to $S$, as in Lemma~\ref{lem:smoothsection}; this is excluded by hypothesis. If both are nonzero, the vanishing of the derivatives instead gives
\[
x\in S\cap B,\qquad dF_x=-u^2\ell_0-v^2\ell_1.
\]
The second equality says that the tangent plane to $S$ at $[x]$ contains $B$, again contradicting transversality. Hence no singular point can occur.
\end{proof}

Conversely, the diagonal form makes smoothness of the two coordinate Eckardt sections necessary. Indeed, a singular point of $\{F+u^2\ell_0=0\}$ remains singular on the fourfold at $v=0$, because the additional derivative $2v\ell_1(x)$ vanishes there. The same argument with $u$ and $v$ interchanged applies to $\{F+v^2\ell_1=0\}$.

For~\eqref{eq:ranktwo}, the parameter map is
\begin{equation}\label{eq:doublecover}
\PP^1\longrightarrow\Lambda,\qquad [u:v]\longmapsto[u^2\ell_0+v^2\ell_1].
\end{equation}
Its branch values are $H_0,H_1$. Varying them leaves the image equal to $\Lambda$, which leads to the following failure of reconstruction.

\Needspace{6\baselineskip}
\begin{theoremB}
There is a two-dimensional algebraic family of pairwise equivariantly nonisomorphic smooth cubic fourfolds with symplectic involutions. Its members have $\rk\widehat q=2$, the same fixed cubic surface $S$ with $\Aut(S)=1$, and the same entire set of polarized isomorphism classes of anti-invariant factors from smooth Eckardt sections.
\end{theoremB}
\begin{proof}
To write the family without ordering its branch values, choose an affine coordinate $r$ on $\Lambda$, so that the plane with coordinate $r$ has equation $a+rb=0$. We use the midpoint $m$ of the branch pair and the square $D$ of its half-difference, and put
\[
q_{m,D}=2uv\,a+(u^2+2muv+Dv^2)b.
\]
For $D\ne0$, the two binary quadratics have no common zero, and the coefficient map has rank two. To identify its branch pair, write the map to $\Lambda$ in the coordinate $t=u/v$ as $r=m+(t+D/t)/2$. The fibre over a finite value of $r$ is then given by
\[
t^2+2(m-r)t+D=0,
\]
whose discriminant is $4((m-r)^2-D)$. Therefore the branch values are $m\pm\sqrt D$, as intended; the map is unramified over infinity.

We next check smoothness by relating this equation to the diagonal form of Lemma~\ref{lem:ranktwosmooth}. Fix $D\ne0$, choose $\eta^2=D$, and change coordinates by
\[
s=u+\eta v,\qquad w=u-\eta v,\qquad
u=\frac{s+w}{2},\qquad v=\frac{s-w}{2\eta}.
\]
With
\[
\ell_+=\frac{a+(\eta+m)b}{2\eta},\qquad
\ell_-=\frac{-a+(\eta-m)b}{2\eta},
\]
substitution gives
\begin{equation}\label{eq:ranktwodiagonalization}
q_{m,D}=s^2\ell_++w^2\ell_-.
\end{equation}
The planes defined by $\ell_+$ and $\ell_-$ have coordinates $m+\eta$ and $m-\eta$, so they are exactly the branch planes. Although this diagonalization uses a square root for each fibre, the equation $F+q_{m,D}=0$ is polynomial in $m,D$ and defines the family without such a choice.

It remains to require that the branch planes lie outside $\Delta$. Accordingly, let $U$ be the complement in the $(m,D)$-plane of the zero locus of
\begin{equation}\label{eq: nonzero pol}
D\prod_{\delta\in\Delta\cap\mathbb A^1}\bigl((m-\delta)^2-D\bigr).
\end{equation}
The polynomial is nonzero, so $U$ is a nonempty open subset of dimension two. By its definition, the two branch planes are distinct and transverse to $S$ for every $(m,D)\in U$. Their common line is the fixed line $B$, which is also transverse to $S$. Therefore equation~\eqref{eq:ranktwodiagonalization} satisfies all the hypotheses of Lemma~\ref{lem:ranktwosmooth}, proving smoothness throughout the family.

To compare the polarized factor classes, we determine which surface--plane pairs arise from smooth sections. Since $c_{m,D}:\PP^1\to\Lambda$ is surjective, equation~\eqref{eq:smoothimage} gives
\[
U_{m,D}=c_{m,D}^{-1}(\Lambda\setminus\Delta),\qquad
c_{m,D}(U_{m,D})=\Lambda\setminus\Delta.
\]
Thus the smooth parameter sets $U_{m,D}$ may vary, but the planes they parametrize are always exactly those in $\Lambda\setminus\Delta$. In particular, the branch values are included, since their sections are smooth. As the polarized isomorphism class of a factor depends only on its surface--plane pair, every member has the same unparametrized Prym image
\[
\{\mathcal P(S,H):H\in\Lambda\setminus\Delta\},
\]
where $\mathcal P(S,H)$ denotes the polarized anti-invariant factor associated with $(S,H)$.

Finally, suppose that two members of the family are equivariantly isomorphic. Such an isomorphism extends to a projectivity preserving the eigenspaces. Since $\Aut(S)=1$, it restricts to the identity on $S$, hence also on $\PP(V)$ and $\Lambda$. The two double covers can therefore differ only by a projectivity of the source, which leaves their unordered branch pairs unchanged. The sum of the two branch values is $2m$, while the square of their difference is $4D$; equality of the branch pairs consequently gives equality of $(m,D)$. Distinct points of $U$ thus yield equivariantly nonisomorphic fourfolds.
\end{proof}

\Needspace{12\baselineskip}
\section{Torsors and surfaces with automorphisms}
\label{sec:extensions}

The reconstruction argument separates the recovery of each polarized factor from the identification of the resulting surfaces. We now allow torsor data in place of the factors, and then drop the condition $\Aut(S)=1$. In the latter case, more factors compensate for the finite ambiguity in identifying the surfaces.
\sectionintrospace

\subsection{Recovering the factors from torsor data}
A Jacobian torsor has no preferred origin, but its translation abelian variety is intrinsic. An affine involution on the torsor therefore still determines a linear involution, from which we can recover the anti-invariant abelian subvariety. To recover its polarization as well, we must know the multiplicity of the supplied numerical class.

\Needspace{10\baselineskip}
\begin{proposition}\label{prop:torsors}
Let $(Y,\iota)$ be a smooth cubic fourfold of the form~\eqref{eq:normalform}, with $\rk\widehat q=3$ and $\Aut(S)=1$, and choose five distinct smooth Eckardt parameters $t_i$. For each $i$, let $P_i$ be a degree-$d_i$ intermediate Jacobian torsor equipped with its induced affine involution and with the numerical polarization class $m_i\theta_i$, where the positive integer $m_i$ is specified. These data determine each polarized factor $(A_{t_i},\lambda_{t_i})$, and hence determine $(Y,\iota)$ up to projective equivariant isomorphism.
\end{proposition}
\begin{proof}
The translation group of $P_i$ is its Albanese variety $J(T_{t_i})$, and the linear part of the given affine involution acts on this group independently of a choice of origin. The numerical polarization class is independent of that choice as well, because translations act trivially on the N\'eron--Severi group.

Since $\NS(J(T_{t_i}))$ is torsion-free, the known multiplicity $m_i$ determines $\theta_i$ from the supplied class. Taking the image of one minus the linear involution and restricting $\theta_i$ to it therefore gives $(A_{t_i},\lambda_{t_i})$. The five recovered factors determine $(Y,\iota)$ by Theorem~\ref{thm:five}.
\end{proof}

For twisted intermediate Jacobians, the required multiplicity is provided by the construction. Voisin constructs the degree-one torsor and proves that its class has order dividing three \cite[Proposition~2.1 and Lemma~2.3]{VoisinTw}. Moreover, the corrected global theta class on a twisted LSV fibre restricts to three times the principal theta class \cite[Lemma~2.1 and Proposition~2.2]{MOc}. Thus these fibres supply precisely the data of Proposition~\ref{prop:torsors}, with the known multiplicity $m_i=3$.

\subsection{A finite bound with automorphisms}

If $S$ has nontrivial automorphisms, Theorem~\ref{thm:cmz} still recovers each surface--plane pair, but there are several identifications with a fixed representative. Since there are only finitely many, sufficiently many common factor classes force one identification to match five points of the conics.

\begin{proposition}\label{prop:finitegroup}
Let $(Y,\iota)$ and $(Y',\iota')$ be smooth cubic fourfolds with symplectic involutions, written in the form~\eqref{eq:normalform}, and assume $\rk\widehat q=\rk\widehat q'=3$. Suppose that their fixed cubic surfaces $S$ and $S'$ are isomorphic, and put $G=\Aut(S)$.

If their smooth Eckardt pencils have more than $4|G|$ polarized anti-invariant factor classes in common, then there is a projective isomorphism $f:Y\xrightarrow{\sim}Y'$ such that $f\circ\iota=\iota'\circ f$. In particular, equality of their entire sets of factor classes implies an equivariant isomorphism.
\end{proposition}
\begin{proof}
Fix an identification of $S'$ with $S$, and write $C=C_q$ and $C'=C_{q'}$ in the resulting common dual space. Theorem~\ref{thm:cmz} associates to each common factor class points $p\in C$, $p'\in C'$ and an element $\gamma\in G$ such that $p=\gamma p'$. Since distinct classes give distinct points $p$, more than $4|G|$ classes force at least five distinct points to use the same $\gamma$. Thus $C\cap\gamma C'$ contains five distinct points.

This forces the two conics to coincide. Indeed, if they spanned different planes, their common points would lie on the intersection line, which meets a smooth conic in at most two points. If they were distinct conics in the same plane, B\'ezout's theorem would instead give at most four common points. Therefore $C=\gamma C'$, and Proposition~\ref{prop:surfaceconic} yields the required equivariant isomorphism.

For the last assertion, observe that the conic has infinitely many points whereas every $G$-orbit is finite. The full Prym image consequently contains infinitely many classes, so equality of the two full images implies the stated finite bound.
\end{proof}

When $G=1$, this argument recovers the sharp bound of Theorem~\ref{thm:five}.

\Needspace{12\baselineskip}
\section{Jacobian compactifications}
\label{sec:lsv}

Once the equivariant cubic is known, it determines the full Jacobian family and its induced involution. We relate this observation to reconstruction from an OG10 fibration and to the moduli construction of cubic pairs. We then prove that the Eckardt-factor family has a Lagrangian closure, before considering maps of fibrations induced by the cubic.
\sectionintrospace

\subsection{Reconstructing the cubic from the full fibration}

Let $B=\PP(V^*\oplus W^*)$ parametrize hyperplanes $K\subset\PP(V\oplus W)$, and let $B^\circ$ be the open locus where $Y\cap K$ is smooth. The intermediate Jacobians form a principally polarized abelian scheme
\[
\pi^\circ:\mathcal J_Y^\circ\longrightarrow B^\circ,
\qquad (\pi^\circ)^{-1}(K)=J(Y\cap K),
\]
with five-dimensional fibres and its natural zero section. To retain information over the singular hyperplane sections, we use a smooth projective hyper-K\"ahler compactification. Such a compactification was constructed for general $Y$ by Laza--Sacc\`a--Voisin \cite{LSV}, and its existence for every smooth cubic fourfold follows from Sacc\`a's theorem \cite[Theorem~1.6]{Sacca}. We fix an associated OG10 Lagrangian fibration $\pi:M_Y\to B$ in the sense of \cite[Definition~4.5]{DMS}; in particular, it compactifies the standard relative intermediate Jacobian over hyperplane sections with at most one ordinary node.

The fibration determines the isotropic class $\eta=c_1(\pi^*\mathcal O_B(1))$ in $H^2(M_Y,\ZZ)$. From this class, the theorem of Dutta--Mattei--Shinder recovers the primitive Hodge lattice of the cubic through an integral Hodge isometry \cite[Theorem~4.6]{DMS}
\[
H^4(Y,\ZZ)_{\mathrm{prim}}(1)\simeq
\eta^\perp/\ZZ\eta,
\qquad \eta^\perp\subset H^2(M_Y,\ZZ).
\]
On the right, the quotient carries the Beauville--Bogomolov form. On the left, the primitive part is orthogonal to the square of the hyperplane class; the Tate twist shifts bidegrees by $(-1,-1)$, and the intersection form is taken with the opposite sign. Since the right-hand side is determined by $\pi$, the isometry recovers the primitive Hodge lattice of $Y$. The Torelli theorem for cubic fourfolds therefore recovers $Y$ itself up to projective isomorphism.

Under the hypotheses $\rk\widehat q=3$ and $\Aut(S)=1$, Theorem~\ref{thm:reconstruction} instead recovers the equivariant cubic from finitely many abstract polarized abelian fourfolds. Once $(Y,\iota)$ is known, its relative Jacobian construction determines $\mathcal J_Y^\circ\to B^\circ$ and the induced involution. Thus the finite data recover the full smooth family, up to the projective change of base induced by the cubic isomorphism.

\subsection{A comparison with moduli of cubic pairs}
\label{sec:moduli-comparison}
The construction of Laza--Pearlstein--Zhang gives a related reconstruction with a moduli interpretation. A cubic with a transverse hyperplane section determines an Eckardt cubic of one higher dimension \cite[Section~1]{LPZ}; for a cubic threefold and a hyperplane, its equation is
\[
Y'=\{G(x_0,\ldots,x_4)+z^2L(x_0,\ldots,x_4)=0\},
\qquad \sigma[x:z]=[x:-z],
\]
where $G$ is cubic and $L$ is linear, and $\sigma$ is anti-symplectic. This construction gives a period map, defined by the fourth cohomology of $Y'$, to a fourteen-dimensional arithmetic quotient of a type IV domain. By \cite[Theorems~3.5 and~5.5]{LPZ}, the period map is an isomorphism onto its image and extends to an isomorphism from a GIT compactification to the Baily--Borel compactification.

The same fourfold also yields polarized anti-invariant factors: for smooth $Y'$, a general invariant hyperplane $h(x)=0$ cuts out an Eckardt threefold. The Torelli theorem for cubic surface pairs \cite[Theorem~6.1]{CMZ} therefore recovers
\[
S_h=\{G=h=0\},\qquad \Pi_h=\{L=h=0\},
\]
as a surface--plane pair in $\{h=0\}\subset\PP^4_x$ from the corresponding factor. Here, however, the surfaces $S_h$ vary among the hyperplane sections of $\{G=0\}$. In the Eckardt pencil of Theorem~\ref{thm:reconstruction}, all the recovered surface--plane pairs instead have the same surface $S$.

\subsection{The induced involution}

We return to the symplectic involution $\iota$ on $Y$. Its action on hyperplanes is
\[
\iota_B[\alpha:\beta]=[\alpha:-\beta],
\qquad \alpha\in V^*,\quad\beta\in W^*.
\]
Since $\iota$ identifies the section over $K$ with the section over $\iota_B(K)$, it induces isomorphisms $J(Y\cap K)\to J(Y\cap\iota_B(K))$. These assemble into an automorphism $f^\circ$ of the smooth Jacobian family satisfying
\[
\pi^\circ\circ f^\circ=\iota_B\circ\pi^\circ.
\]
In particular, a hyperplane fixed by $\iota_B$ gives an involution of a single Jacobian. This applies to every smooth Eckardt hyperplane $K_t=\PP(V\oplus\CC t)$: for $t\in U_q$, the restriction of $f^\circ$ to $J(T_t)$ is induced by $\iota_t$, and its anti-invariant factor is $A_t$.

Because the smooth Jacobian family is dense in $M_Y$, this automorphism extends uniquely to a birational involution $f:M_Y\dashrightarrow M_Y$ \cite[Section~3.1]{Sacca}. Moreover, any two associated compactifications contain the same smooth family with the same action, and are consequently equivariantly birational. It follows that the finite data of Theorem~\ref{thm:reconstruction} determine the equivariant birational class of $(M_Y,f)$, up to the projective change of base induced by the cubic isomorphism.

Moreover, $f$ is symplectic because $\iota$ acts trivially on $H^{3,1}(Y)$ \cite[Theorem~5.1 and its proof]{MM}, and it is nontrivial because its action on $B$ is nontrivial. On the other hand, the symplectic rigidity theorem of Giovenzana--Grossi--Onorati--Veniani states that every finite-order symplectic regular automorphism of a smooth OG10 manifold is the identity \cite[Theorem~1.1]{GGOV}. Therefore, although $M_Y$ is smooth, $f$ cannot be regular on it or on any smooth hyper-K\"ahler birational model.

\subsection{The closure of the Eckardt-factor family}
Identify the Eckardt line in $B$ with $\PP(W)$ by $t\mapsto K_t$, and write $U=U_q$. The factors form an abelian subscheme $\mathcal A_U\to U$, closed in $\mathcal J_Y^\circ|_U$, of total dimension $1+4=5$. To show that it is Lagrangian, we will compose $f^\circ$ with fibrewise inversion and show that the resulting anti-symplectic involution fixes it pointwise.

\Needspace{8\baselineskip}
\begin{proposition}\label{prop:prymclosure}
Let $M_Y$ be an associated OG10 compactification of the intermediate Jacobian family of $Y$. The scheme-theoretic closure
\[
Z_Y=\overline{\mathcal A_U}\subset M_Y
\]
is an integral projective Lagrangian subvariety.
\end{proposition}
\begin{proof}
Let $\omega$ be a holomorphic symplectic form on $M_Y$, and let $\nu$ be fibrewise inversion on $\mathcal J_Y^\circ$. By \cite[Section~3.1]{Sacca}, inversion extends birationally to $M_Y$. Since $H^0(M_Y,\Omega^2_{M_Y})=\CC\omega$, it follows that $\nu^*\omega=c\omega$ for some $c\in\CC^*$.

To determine this scalar, consider a point $p$ of the zero section. The section splits the tangent space as $T_p\mathcal J_Y^\circ=H_p\oplus V_p$, where $H_p$ is tangent to the section and $V_p$ to the fibre. Inversion fixes the zero section and acts by multiplication by $-1$ along the fibres, so $d\nu_p$ is $+\mathrm{id}$ on $H_p$ and $-\mathrm{id}$ on $V_p$. Since $V_p$ is Lagrangian and $\omega_p$ is nondegenerate, there exist $h\in H_p$ and $v\in V_p$ with $\omega_p(h,v)\ne0$. Evaluating the two forms on this pair gives
\[
c\,\omega_p(h,v)
=(\nu^*\omega)_p(h,v)
=\omega_p(h,-v)
=-\omega_p(h,v),
\]
and hence $\nu^*\omega=-\omega$.

We can now combine inversion with the induced involution. The map $f^\circ$ is symplectic \cite[Theorem~5.1 and its proof]{MM} and acts by group isomorphisms between fibres, so it commutes with $\nu$. Consequently, $g=\nu\circ f^\circ$ is a regular anti-symplectic involution of $\mathcal J_Y^\circ$. On $J(T_t)$, for $t\in U$, it acts as $-\iota_t$. Moreover, the identity $\iota_t(1-\iota_t)=-(1-\iota_t)$ shows that $\iota_t$ acts by multiplication by $-1$ on $A_t=\im(1-\iota_t)$. Thus $g$ fixes the whole family $\mathcal A_U$ pointwise.

The local fixed-locus argument from the proof of \cite[Lemma~5.4]{BCS} now gives the required vanishing. Indeed, for tangent vectors $u,v$ to $\mathcal A_U$, the differential of $g$ is the identity, so
\[
\omega(u,v)=\omega(dg(u),dg(v))
=(g^*\omega)(u,v)=-\omega(u,v).
\]
It follows that $\omega$ vanishes on $\mathcal A_U$. As the family has dimension $1+4=5=\frac12\dim M_Y$, it is Lagrangian.

Finally, we pass to its closure. The base $U$ is an irreducible smooth curve, so the abelian scheme $\mathcal A_U$ is integral. Its scheme-theoretic closure $Z_Y$ is therefore an integral projective subvariety of dimension five, containing $\mathcal A_U$ as a dense smooth open subset. Since the restriction of $\omega$ to $Z_Y^{\mathrm{reg}}$ vanishes on this dense open subset, it vanishes everywhere on $Z_Y^{\mathrm{reg}}$. Thus $Z_Y$ is Lagrangian as claimed.
\end{proof}

For any associated compactification $M_Y$, the closure contains the same dense open family $\mathcal A_U$. Its birational class over $\PP(W)$ is therefore independent of the compactification.

\subsection{Maps induced by the cubic}
Reconstructing the cubic leaves a further question: is a given polarized isomorphism of Jacobian fibrations, preserving their zero sections, induced by a projectivity of the cubics? Fibrewise inversion shows that this stronger conclusion need not hold.

Indeed, $\nu$ covers the identity on $B^\circ$, preserves the zero section and the principal polarizations, and commutes with $f^\circ$. If it were induced by a projectivity $g$ of $Y$, the dual projectivity would fix every point of $B^\circ$, hence would be the identity on $B$. This would force $g$ to be the identity on the ambient projective space and therefore to induce the identity on every Jacobian, contrary to the nontriviality of $\nu$. On the other hand, the proof of Proposition~\ref{prop:prymclosure} shows that inversion is anti-symplectic.

If the zero section need not be preserved, translations provide another possible source of birational maps. Namely, whenever there is a nonzero rational section $s:B\dashrightarrow M_Y$, translation by $s$ acts on each smooth fibre where it is defined by $a\mapsto a+s(K)$ and preserves the numerical polarization class. Since $f$ acts by group isomorphisms on the smooth fibres, this translation commutes with $f$ exactly when
\[
f\circ s=s\circ\iota_B
\]
as rational maps. The availability of such translations is governed by the group of rational sections, which Sacc\`a's theorem identifies with $H^{2,2}(Y,\ZZ)_{\mathrm{prim}}$ \cite[Theorem~5.1]{Sacca}.

\begingroup
\raggedright
\bibliographystyle{alphaurl}
\bibliography{references}
\endgroup
\end{document}